\documentclass[12pt]{amsart}
\calclayout

\usepackage{amsmath}
\usepackage{amssymb}
\usepackage{amsfonts}
\usepackage{amsthm}
\usepackage{mathrsfs}
\usepackage[all]{xy}
\usepackage{enumerate}
\usepackage{graphicx}
\usepackage{romannum}
\usepackage{xcolor}

\usepackage{physics}
\usepackage{relsize}
\usepackage{hyperref}

\usepackage[autostyle]{csquotes}

\def\RR{\mathbb{R}} 
\def\CC{\mathbb{C}} 
\def\NN{\mathbb{N}} 
 
\def\dist{\mathrm{dist}} 
 
\def\O{\Omega} 
\def\Levi{\mathscr{L}} 

\newtheorem{thm}{Theorem}[section]
\newtheorem{lem}[thm]{Lemma}

\newtheorem{cor}[thm]{Corollary}

\theoremstyle{definition}
\newtheorem{defn}[thm]{Definition}

\theoremstyle{remark}
\newtheorem{rmk}[thm]{Remark}

\numberwithin{equation}{section}

\def\BigRoman{\uppercase\expandafter{\romannumeral\number\count 255 }}
\def\Romannumeral{\afterassignment\BigRoman\count255=}

\begin{document}
	
	\renewcommand{\thepage}{\arabic{page}}
	
	\title{On the Steinness index}         
	\author{Jihun Yum}        
	\address{Department of Mathematics, Pohang University of Science and Technology, Pohang, 790-784, Republic of Korea}
	\email{wadragon@postech.ac.kr}
	\thanks{This research was supported by the SRC-GAIA (NRF-2011-0030044) through the National Research Foundation of Korea (NRF) funded by the Ministry of Education.}

	\begin{abstract}
		We introduce the concept of Steinness index related to the Stein neighborhood basis. We then show several results: (1) The existence of Steinness index is equivalent to that of strong Stein neighborhood basis. (2) On the Diederich-Forn{\ae}ss worm domains in particular, we present an explicit formula relating the Steinness index to the well-known Diederich-Forn{\ae}ss index. (3) The Steinness index is 1 if a smoothly bounded pseudoconvex domain admits finitely many boundary points of infinite type.
	\end{abstract}
	
	\maketitle

	\section{\bf Introduction} \label{section introduction}
	
	Let $\O \subset \CC^n$$(n \ge 2)$ be a bounded domain with smooth boundary. 
	A smooth function $\rho$ defined on a neighborhood $V$ of $\overline{\O}$ is called a (global) {\it defining function} of $\O$ if $\O = \{ z \in V : \rho(z)<0 \}$ and $d\rho(z) \neq 0$ for all $z \in \partial \O$.

	\subsection{Diederich-Forn{\ae}ss index}
	The {\it Diederich-Forn{\ae}ss exponent of $\rho$} is defined by 
	$$\eta_{\rho} := \sup\{ \eta \in (0,1) : -(-\rho)^{\eta} \text{ is strictly plurisubharmonic on } \O \}.$$
	If there is no such $\eta$, then we define $\eta_{\rho} = 0$. 
	{\it The Diederich-Forn{\ae}ss index of $\O$} is defined by 
	$$DF(\O) := \sup \text{ } \eta_{\rho},$$ 
	where the supremum is taken over all defining functions $\rho$. We say that the Diederich-Forn{\ae}ss index of $\O$ exists if $DF(\O) \in (0, 1]$. If $DF(\O)$ exists, then there exists a bounded strictly plurisubharmonic exhaustion function on $\O$. In other words, $\O$ becomes a hyperconvex domain. In 1977, Diederich and Forn{\ae}ss (\cite{DieForn1977}) proved that $C^2$-smoothness of $\partial \O$ implies the existence of $DF(\O)$.

	\subsection{Steinness index}
	
	We introduce the following definition. The {\it Steinness exponent of $\rho$} is defined by 
	\begin{align*}
		\widetilde{\eta}_{\rho} := \inf\{ \widetilde{\eta} > 1 : \rho^{\widetilde{\eta}} \text{ is strictly plurisubharmonic on } \overline{\Omega}^{\complement} \cap U  &\\
		\text{ for some neighborhood } U \text{ of } & \partial \O  \},
	\end{align*}
	where $\overline{\O}^{\complement}:= \CC^n \setminus \overline{\O}$.
	If there is no such $\widetilde{\eta}$, then we define $\widetilde{\eta}_{\rho} = \infty$. 
	The {\it Steinness index of $\Omega$} is defined by 
	$$S(\Omega) := \inf \text{ } \widetilde{\eta}_{\rho},$$ 
	where the infimum is taken over all defining functions $\rho$.  We say that the Steinness index of $\O$ exists if $S(\O) \in [1,\infty)$.
	$\overline{\O}$ is said to have a {\it Stein neighborhood basis} if for any neighborhood $V_1$ of $\overline{\O}$, there exists a pseudoconvex domain $V_2$ such that $\overline{\O} \subset V_2 \subset V_1$.
	If $S(\O)$ exists, then there exist a defining function $\rho$ and $\eta_2 \in (1,\infty)$ such that $\rho^{\eta_2}$ is strictly plurisubharmonic on $\overline{\O}^{\complement} \cap U$. Thus $\overline{\O}$ has a Stein neighborhood basis. In contrast with Diederich-Forn{\ae}ss index, the smoothness of boundary does not implies the existence of $S(\O)$; Diederich-Forn{\ae}ss worm domains provide an example (\cite{DieForn1977-2}). 
	In section \ref{section equivalent definition}, we characterize the Steinness index by means of a differential inequality on the set of all weakly pseudoconvex boundary points (Theorem \ref{S main thm}). This theorem plays a crucial role in this paper.

	\subsection{Strong Stein neighborhood basis}
	
	
	$\overline{\O} \subset \subset \CC^n$ is said to have a {\it strong Stein neighborhood basis} if there exist a defining function $\rho$ of $\O$ and $\epsilon_0 > 0$ such that 
	\begin{equation*}
	\O_{\epsilon} := \{ z \in \CC^n : \rho(z) < \epsilon \}
	\end{equation*}
	is pseudoconvex for all $0 \le \epsilon < \epsilon_0$. This implies the existence of a Stein neighborhood basis. In section \ref{section strong Stein neighborhood basis}, we show that the existence of $S(\O)$ is actually equivalent to the existence of a strong Stein neighborhood basis (Theorem \ref{Sindex eqi SSnbhd}).

	\subsection{Worm domains}
	
	In 1977, Diederich and Forn{\ae}ss (\cite{DieForn1977-2}) constructed bounded smooth domains in $\CC^2$ whose Diederich-Forn{\ae}ss indices are strictly less than one. These examples are called worm domains and the only known domains in $\CC^n$ which have non-trivial Diederich-Forn{\ae}ss indices.
	Therefore, worm domains are worth to be studied.
	
	Recently, Liu calculated the exact value of the Diederich-Forn{\ae}ss index of worm domains (Definition \ref{worm defn}) in 2017. In section \ref{section Steinness index of worm domains}, exploiting the idea of \cite{Liu2017}, we obtain a calculation of the exact values of the Steinness index of worm domains (Theorem \ref{S worm}). More precisely, the result is as follows.

	\begin{thm} \label{main}
		If $\O_{\beta}$ $(\beta > \frac{\pi}{2})$ is a worm domain, then the following 4 conditions are equivalent:
		\begin{itemize}
			\item[1.] $\frac{1}{2} < DF(\O_{\beta}) < 1$.
			\item[2.] $\O_{\beta}$ admits the Steinness index.
			\item[3.] $\O_{\beta}$ admits the Stein neighborhood basis.
			\item[4.] $\O_{\beta}$ admits a strong Stein neighborhood basis.	
		\end{itemize}
		Moreover, if one of the above conditions holds, then 
		$$ \frac{1}{DF(\O_{\beta})} + \frac{1}{S(\O_{\beta})} = 2. $$
	\end{thm}

	\subsection{Sufficient conditions for ${\bf DF(\O)=1}$ and ${\bf S(\O)=1}$}

	If a bounded domain $\O \subset \CC^n$ is strongly pseudoconvex, then there exists a strictly plurisubharmonic defining function of $\O$, which implies $DF(\O)=1$ and $S(\O)=1$.

	In fact, more is known: If a smoothly bounded pseudoconvex domain $\O$ is B-regular (i.e., for every $p \in \partial \O$, there exists a continuous peak function at $p$), then $DF(\O)=1$ and $S(\O)=1$ (\cite{Sibony1987}, \cite{Sibony1991}).
	Since finite type domains in the sense of D'Angelo are B-regular by Catlin (\cite{Catlin1984}), both indices are equal to 1. In section \ref{section Steinness index of finite type domains}, we give an alternative proof of this fact using modified Theorem \ref{S main thm} (Corollary \ref{finite type}). Moreover, if the set of all infinite type boundary points is finite, then the Steinness index is 1 (Corollary \ref{finite infinite type points}). In section \ref{section Steinness index of convex domains}, we also demonstrate that if $\O \subset\subset \CC^n$ is a $C^1$-smooth convex domain, then $DF(\O)=1$ and $S(\O)=1$ (Corollary \ref{convex domain index 1}). 	
	For further results, we refer the reader to \cite{DieForn1977-3}, \cite{ForHer2007}, \cite{ForHer2008}, \cite{Harring2018} and \cite{KranLiuPelo2016}.

\vspace{5mm}


	\section{\bf Preliminary}

	We first fix the notation of this paper, unless otherwise mentioned.
	\begin{itemize}
		\item[$\bullet$] $\O$ : a bounded pseudoconvex domain with smooth boundary in $\CC^n$.
		\item[$\bullet$] $\rho$ : a defining function of $\O$.
		\item[$\bullet$] $\Sigma$ : the set of all weakly pseudoconvex points in $\partial \O$.
		\item[$\bullet$] $\Sigma_{\infty}$ : the set of all infinite type points in $\partial \O$.
		\item[$\bullet$] $g$ : the standard Euclidean complex Hermitian metric in $\CC^n$.
		\item[$\bullet$] $\nabla$ : the Levi-Civita connection of $g$.
		\item[$\bullet$] $\nabla \rho$ : the real gradient of $\rho$.
		\item[$\bullet$] $U$ : a tubular neighborhood of $\partial \O$.
		\item[$\bullet$] $\Levi_{\rho}$ : the Levi-form of $\rho$.
	\end{itemize}
	Define
	\begin{equation*}
		N_{\rho} := \frac{1}{\sqrt{ \sum^n_{j=1} |\frac{\partial \rho}{\partial z_j}|^2}} \sum^n_{j=1} \frac{\partial \rho}{\partial \bar{z}_j} \frac{\partial}{\partial z_j}.
	\end{equation*}
	Let $J$ be the complex structure of $\CC^n$ and $T_p (\partial \O)$ be the real tangent space of $\partial \O$ at $p \in \partial \O$. Let $T_p^c(\partial \O) := J(T_p (\partial \O)) \cap T_p (\partial \O)$. Then the complexified tangent space of $T_p^c(\partial \O)$, $\CC T_p^c(\partial \O) := \CC \otimes T_p^c(\partial \O)$, can be decomposed into the holomorphic tangent space $T^{1,0}_p(\partial \O)$ and the anti-holomorphic tangent space $T^{0,1}_p(\partial \O)$. We call $X$ a $(1,0)$ tangent vector if $X \in T^{1,0}_p(\partial \O)$.
	
	Let $X,Y,Z$ be complex vector fields in $\CC^n$. 
	A direct calculation implies the following properties.
	\begin{align*}
		g(N_{\rho}, N_{\rho}) = \frac{1}{2}, \phantom{aa}
		N_{\rho} \rho = \sqrt{\sum^n_{j=1} \left| \frac{\partial \rho}{\partial z_j} \right|^2} = \frac{\lVert \nabla \rho \rVert}{2},  
	\end{align*}
	$$ N_{\rho} + \overline{N}_{\rho} = 2 \Re(N_{\rho}) = \frac{\nabla \rho}{\lVert \nabla \rho \rVert},   $$
	$$ \Levi_{\rho}(X,Y) = g(\nabla_X \nabla \rho, Y) = X(\overline{Y} \rho) - (\nabla_X \overline{Y}) \rho,  $$
	\begin{align*}
		Z g(X,Y) = g(\nabla_Z X, Y) + g(X, \nabla_{\overline{Z}}Y), \phantom{aa}
		Z \rho = g(\nabla \rho, \overline{Z}),
	\end{align*}
	$$ \nabla_X\nabla_Y - \nabla_Y\nabla_X - \nabla_{[X,Y]} = 0. $$
	For $p \in \partial \O$ whenever we mention $\lim\limits_{z \rightarrow p}$, it means $z$ approaches $p$ along the real normal direction, and \enquote{smooth} means $C^{\infty}$-smooth, although $C^3$-smoothness suffices for our purpose.

	Now, we introduce two lemmas which we will use in the next section.

	\begin{lem} [\cite{Liu2017}] \label{orthgonal}
		Let $\Omega \subset\subset \CC^n$ be a pseudoconvex domain with smooth boundary, and $\rho$ be a defining function of $\Omega$. Suppose that $L$ is a $(1,0)$ tangent vector field so that $\Levi_{\rho}(L,L) = 0$ at $p \in \partial \O$. Assume that $T_j$ $(1 \le j \le n-2)$ are $(1,0)$ tangent vector fields and $T_1, T_2, \cdots, T_{n-2}, L$ are orthogonal at $p$. Then $\Levi_{\rho}(L, T_j) = 0$ for $1 \le j \le n-2$ at $p$.
	\end{lem}

	\begin{lem} \label{keylem}
		Let $\Omega \subset\subset \CC^n$ be a pseudoconvex domain with smooth boundary, and $\rho$ be a defining function of $\Omega$. Let $\psi$ be a smooth function defined on $U$ and $\widetilde{\rho} = \rho e^{\psi}$. We denote $\widetilde{N} = N_{\widetilde{\rho}}$, $N = N_{\rho}$. Let $\widetilde{L}$, $L$ be $(1,0)$ tangent vector fields on $U$ such that $\widetilde{L} = L$ on $\partial \O$, and  $\widetilde{L} \widetilde{\rho} = 0$, $L\rho = 0$ on $U$. 
		Define
		$$ \Sigma_L := \{ p \in \partial \O : \Levi_{\rho}(L,L)(p) = 0 \}. $$
		Then
		$$ \frac{\Levi_{\widetilde{\rho}}(\widetilde{L}, \widetilde{N})}{\lVert \nabla \widetilde{\rho} \rVert}  = \frac{\Levi_{\rho}(L,N)}{\lVert \nabla \rho \rVert} + \frac{1}{2}(L\psi) $$
		and
		\begin{equation} \label{2.111}
		\frac{\widetilde{N} \Levi_{\widetilde{\rho}}(\widetilde{L},\widetilde{L})}{\lVert \nabla \widetilde{\rho} \rVert} = \frac{N\Levi_{\rho}(L,L)}{\lVert \nabla \rho \rVert} + \frac{1}{2}\Levi_{\psi}(L,L) - \frac{1}{2}|L\psi|^2 - 2 Re \left[ \frac{\Levi_{\rho}(L,N)}{\lVert \nabla \rho \rVert}(\overline{L}\psi)  \right]
		\end{equation}
		on $\Sigma_L$. In particular, 
		$$ \frac{|\Levi_{\widetilde{\rho}}(\widetilde{L}, \widetilde{N})|^2}{\lVert \nabla \widetilde{\rho} \rVert^2} + \frac{1}{2} \frac{\widetilde{N} \Levi_{\widetilde{\rho}}(\widetilde{L},\widetilde{L})}{\lVert \nabla \widetilde{\rho} \rVert} 
		= \frac{|\Levi_{\rho}(L,N)|^2}{\lVert \nabla \rho \rVert^2} + \frac{1}{2} \frac{N\Levi_{\rho}(L,L)}{\lVert \nabla \rho \rVert} + \frac{1}{4}\Levi_{\psi}(L,L) $$
		on $\Sigma_L$.
	\end{lem}
	\begin{proof}
		Since $\widetilde{L} = L$, $\widetilde{N} = N$, and $\lVert \nabla \widetilde{\rho} \rVert = e^{\psi} \lVert \nabla \rho \rVert$ on $\partial \O$, 
		\begin{align*}
		\Levi_{\widetilde{\rho}}(\widetilde{L}, \widetilde{N}) 
		&= \Levi_{\widetilde{\rho}}(L,N) \\
		&= e^{\psi} \Levi_{\rho}(L,N) + \rho \Levi_{e^{\psi}}(L,N) + (L\rho)(\overline{N}e^{\psi}) + (L e^{\psi})(\overline{N}\rho) \\
		&= e^{\psi} \Levi_{\rho}(L,N) + e^{\psi}(L\psi)\frac{\lVert \nabla \rho \rVert}{2}.
		\end{align*}
		Therefore, we have
		$$ \frac{\Levi_{\widetilde{\rho}}(\widetilde{L}, \widetilde{N})}{\lVert \nabla \widetilde{\rho} \rVert}  = \frac{\Levi_{\rho}(L,N)}{\lVert \nabla \rho \rVert} + \frac{1}{2} (L \psi) $$
		on $\Sigma_L$.
		Now, we prove (\ref{2.111}) holds at $p \in \Sigma_L$.
		\begin{align} \label{2.222}
		\widetilde{N} \Levi_{\widetilde{\rho}}(\widetilde{L},\widetilde{L}) 
		&= N\Levi_{\widetilde{\rho}}(\widetilde{L},\widetilde{L}) 
		= N g(\nabla_{\widetilde{L}} \nabla \widetilde{\rho}, \widetilde{L}) 
		= g(\nabla_N \nabla_{\widetilde{L}} \nabla \widetilde{\rho}, \widetilde{L}) + g(\nabla_{\widetilde{L}} \nabla \widetilde{\rho}, \nabla_{\overline{N}} \widetilde{L}) \\
		&= g(\nabla_{\widetilde{L}} \nabla_N \nabla \widetilde{\rho}, \widetilde{L})  +  g(\nabla_{[N,\widetilde{L}]} \nabla \widetilde{\rho}, \widetilde{L})  +  g(\nabla_{\widetilde{L}} \nabla \widetilde{\rho}, \nabla_{\overline{N}} \widetilde{L}). \nonumber
		\end{align}
		If $L$ vanishes at $p$, then the last three terms of (\ref{2.222}) are all zero and so $\widetilde{N} \Levi_{\widetilde{\rho}}(\widetilde{L},\widetilde{L}) = N\Levi_{\rho}(L,L) = 0$. Thus (\ref{2.111}) holds at $p$. Therefore, we assume that $L \neq 0$ at $p$ and normalize $L$ and $\widetilde{L}$ as $g(L,L) = g(\widetilde{L}, \widetilde{L}) = \frac{1}{2}$. This normalization is for convenience and does not affect the result. Now, we will compute the last three terms of (\ref{2.222}) one by one.
		Let $\{\sqrt{2}\widetilde{T_1}, \cdots, \sqrt{2}\widetilde{T}_{n-2}, \sqrt{2}\widetilde{L} \} $ be an orthonormal basis of $T^{1,0}_p(\partial \O)$. Then first,
		\begin{align*}
		& g(\nabla_{\widetilde{L}} \nabla \widetilde{\rho}, \nabla_{\overline{N}} \widetilde{L}) \\
		= & g \left(  \nabla_{\widetilde{L}} \nabla \widetilde{\rho}, \sum_{j=1}^{n-2}g(\nabla_{\overline{N}} \widetilde{L},\sqrt{2}\widetilde{T_j})\sqrt{2}\widetilde{T_j}  + g(\nabla_{\overline{N}} \widetilde{L}, \sqrt{2}\widetilde{L})\sqrt{2}\widetilde{L}  +  g(\nabla_{\overline{N}} \widetilde{L}, \sqrt{2}N)\sqrt{2}N   \right) \\
		= & \sum_{j=1}^{n-2} 2 \overline{g(\nabla_{\overline{N}} \widetilde{L}, \widetilde{T_j})} \Levi_{\widetilde{\rho}}(\widetilde{L}, \widetilde{T_j}) + 
		2 \overline{g(\nabla_{\overline{N}} \widetilde{L}, \widetilde{L})} \Levi_{\widetilde{\rho}}(\widetilde{L}, \widetilde{L}) +
		2 \overline{g(\nabla_{\overline{N}} \widetilde{L}, N)} \Levi_{\widetilde{\rho}}(\widetilde{L}, N) \\
		= & 2 \overline{g(\nabla_{\overline{N}} \widetilde{L}, N)} \Levi_{\widetilde{\rho}}(\widetilde{L}, N) \\
		= & 2 \overline{g(\nabla_{\overline{N}} \widetilde{L}, N)} \Levi_{\widetilde{\rho}}(L, N).
		\end{align*}
		Here, we used Lemma \ref{orthgonal} in the third equality above. Second, by the same argument as above,
		\begin{align*}
		g(\nabla_{[N,\widetilde{L}]} \nabla \widetilde{\rho}, \widetilde{L}) 
		=  \overline{g(\nabla_{\widetilde{L}} \nabla \widetilde{\rho}, [N,\widetilde{L}])} 
		=  2 g([N,\widetilde{L}], N) \overline{\Levi_{\widetilde{\rho}}(\widetilde{L}, N)}.
		\end{align*}
		Third,
		\begin{align*}
		& g(\nabla_{\widetilde{L}} \nabla_N \nabla \widetilde{\rho}, \widetilde{L}) \\
		= & g(\nabla_{L} \nabla_N \nabla \widetilde{\rho}, L) \\
		= & g(\nabla_{N} \nabla_L \nabla \widetilde{\rho}, L) - g(\nabla_{[N,L]} \nabla \widetilde{\rho}, L) \\
		= & N g(\nabla_L \nabla \widetilde{\rho}, L) - g(\nabla_L \nabla \widetilde{\rho}, \nabla_{\overline{N}} L)  - g(\nabla_{[N,L]} \nabla \widetilde{\rho}, L) \\
		= & N \Levi_{\widetilde{\rho}}(L,L) - 2 \overline{g(\nabla_{\overline{N}} L, N)} \Levi_{\widetilde{\rho}}(L, N) - 2 g([N,L], N) \overline{\Levi_{\widetilde{\rho}}(L, N)}.
		\end{align*}
		Therefore, we have
		\begin{align*}
		\widetilde{N} \Levi_{\widetilde{\rho}}(\widetilde{L},\widetilde{L}) 
		&= N \Levi_{\widetilde{\rho}}(L,L)  \\
		&+ 2 \left( \overline{g(\nabla_{\overline{N}} \widetilde{L}, N)} - \overline{g(\nabla_{\overline{N}} L, N)} \right) \Levi_{\widetilde{\rho}}(L, N)  \\
		&+ 2 \left( g([N,\widetilde{L}], N) - g([N,L], N) \right) \overline{\Levi_{\widetilde{\rho}}(L, N)}.
		\end{align*}
		Now,
		\begin{align}
		\Levi_{\widetilde{\rho}}(L,L) &= e^{\psi} \Levi_{\rho}(L,L) + \rho \Levi_{e^{\psi}}(L,L) + (L\rho)(\overline{L}e^{\psi}) + (L e^{\psi})(\overline{L}\rho), \nonumber \\
		N \Levi_{\widetilde{\rho}}(L,L) &= e^{\psi} \left( N \Levi_{\rho}(L,L) + (N\rho)\Levi_{\psi}(L,L) + (N\rho)|L\psi|^2  \right) \label{1},
		\end{align}
		and
		\begin{align*}
		& \Levi_{\rho}(N,L) = N(\overline{L}\rho) - (\nabla_N \overline{L})\rho = - (\nabla_N \overline{L})\rho \\
		\Rightarrow \phantom{a} & - \Levi_{\rho}(L,N) = (\nabla_{\overline{N}} L)\rho = 2 g(\nabla_{\overline{N}} L, N)(N\rho) = g(\nabla_{\overline{N}} L, N) \lVert \nabla \rho \rVert \\
		\Rightarrow \phantom{a} & g(\nabla_{\overline{N}} L, N) = - \frac{\Levi_{\rho}(L,N)}{\lVert \nabla \rho \rVert}. \\
		\end{align*}
		By the same argument,
		\begin{align}
		& g(\nabla_{\overline{N}} \widetilde{L}, N) = - \frac{\Levi_{\widetilde{\rho}}(\widetilde{L},N)}{\lVert \nabla \widetilde{\rho} \rVert} = - \frac{\Levi_{\rho}(L,N)}{\lVert \nabla \rho \rVert} - \frac{1}{2}(L \psi) \nonumber \\
		\Rightarrow \phantom{a} & g(\nabla_{\overline{N}} \widetilde{L}, N) - g(\nabla_{\overline{N}} L, N) = - \frac{1}{2}(L\psi) \label{2},
		\end{align}
		and
		\begin{align*}
		g([N,L],N) 
		= g(\nabla_N L, N) - g(\nabla_L N, N) 
		= g(\nabla_{\overline{N}} L, N) + g(\nabla_{N-\overline{N}} L, N) -  g(\nabla_L N, N), \\
		g([N,\widetilde{L}],N) 
		= g(\nabla_{N} \widetilde{L}, N) - g(\nabla_{\widetilde{L}} N, N) 
		= g(\nabla_{\overline{N}} \widetilde{L}, N) + g(\nabla_{N-\overline{N}} \widetilde{L}, N) -  g(\nabla_{\widetilde{L}} N, N).
		\end{align*}
		Here, $g(\nabla_L N, N) = g(\nabla_{\widetilde{L}} N, N)$ on $\Sigma_L$ and since $N-\overline{N}$ is a real tangent vector field on $\partial \O$, $g(\nabla_{N-\overline{N}} L, N) = g(\nabla_{N-\overline{N}} \widetilde{L}, N)$ on $\Sigma_L$. Therefore,
		\begin{align} \label{3}
		g([N,\widetilde{L}],N) - g([N,L],N) = g(\nabla_{\overline{N}} \widetilde{L}, N) - g(\nabla_{\overline{N}} L, N) = - \frac{1}{2}(L\psi).
		\end{align}
		Combining (\ref{1}), (\ref{2}), and (\ref{3}), we have
		\begin{align*}
		&\widetilde{N} \Levi_{\widetilde{\rho}}(\widetilde{L},\widetilde{L}) \\
		= & e^{\psi} \left( N \Levi_{\rho}(L,L) + (N\rho)\Levi_{\psi}(L,L) + (N\rho)|L\psi|^2  \right) \\
		&- (\overline{L}\psi) \left( e^{\psi} \Levi_{\rho}(L,N) + e^{\psi}(L\psi)(N\rho)\right) \\
		&- (L\psi) \left( e^{\psi} \overline{\Levi_{\rho}(L,N)} + e^{\psi}(\overline{L}\psi)(N\rho)\right) \\
		= & e^{\psi} \left[ N \Levi_{\rho}(L,L) + (N\rho)\Levi_{\psi}(L,L) - (N\rho)|L\psi|^2 - 2 \Re\left( \Levi_{\rho}(L,N)(\overline{L}\psi) \right) \right].
		\end{align*}
		Since $\lVert \nabla \widetilde{\rho} \rVert = e^{\psi} \lVert \nabla \rho \rVert$ on $\partial \O$, we conclude
		$$ \frac{\widetilde{N} \Levi_{\widetilde{\rho}}(\widetilde{L},\widetilde{L})}{\lVert \nabla \widetilde{\rho} \rVert} 
		= \frac{N\Levi_{\rho}(L,L)}{\lVert \nabla \rho \rVert} + \frac{1}{2}\Levi_{\psi}(L,L) - \frac{1}{2}|L\psi|^2 - 2 \Re \left[ \frac{\Levi_{\rho}(L,N)}{\lVert \nabla \rho \rVert}(\overline{L}\psi)  \right] $$
		at $p \in \Sigma_L$.
	\end{proof}

\vspace{5mm}


	\section{\bf Equivalent Definition for Steinness index} \label{section equivalent definition}

	In order to check whether $S(\O)$ exists, by its definition, it is necessary to find a defining function $\rho$ and $\eta_2 >1$ such that $\rho^{\eta_2}$ is strictly plurisubharmonic on $\overline{\O}^{\complement} \cap U$.
	In this section, we replace this condition on $\overline{\O}^{\complement} \cap U$ to another condition on $\Sigma$, so that we only need to check it at weakly pseudoconvex boundary points. Now, we introduce the main theorem of this section.

	\begin{thm} \label{S main thm}
		Let $\Omega \subset\subset \CC^n$ be a pseudoconvex domain with smooth boundary, and $\rho$ be a defining function of $\Omega$. Let $L$ be an arbitrary $(1,0)$ tangent vector field on $\partial \O$. Define
		$$ \Sigma_L := \{ p \in \partial \O : \Levi_{\rho}(L,L)(p) = 0 \}. $$
		Let $\eta_{\rho} $ be the infimum of $\eta_2 \in (1,\infty)$ satisfying 
		$$ \frac{1}{\eta_2-1} \frac{|\Levi_{\rho}(L,N)|^2}{\lVert \nabla \rho \rVert^2} - \frac{1}{2} \frac{N\Levi_{\rho}(L,L)}{\lVert \nabla \rho \rVert} \le 0 $$
		on $\Sigma_L$ for all $L$. Here, $N = N_{\rho}$ and we extend $L$ so that $L \rho = 0$ on $U$.
		Then 
		$$S(\O) = \inf \text{ } \eta_{\rho}$$
		where the infimum is taken over all smooth defining functions $\rho$.
	\end{thm}

	The following sequence of lemmas are essential towards the proof of the theorem above.

	\begin{lem} \label{S lem1}
		Let $\Omega \subset\subset \CC^n$ be a pseudoconvex domain with smooth boundary, and $\rho$ be a defining function of $\Omega$. For $p \in \partial \O$, let $U_p$ be a neighborhood of $p$ in $\CC^n$. Let $L$ be a smooth $(1,0)$ tangent vector field in $U_p$ such that $L \rho = 0$ and $\Levi_{\rho}(L,L) = 0$. Fix $\eta_2 \in (1,\infty)$. We denote $N_{\rho}$ by $N$.   Then 
		$$ \Levi_{\rho^{\eta_2}}(aL + bN, aL + bN) > 0 $$
		for all $(a,b) \in \CC^2 \setminus (0,0)$ on $\overline{\O}^{\complement} \cap U_p$ implies
		$$ \frac{1}{\eta_2-1} \frac{|\Levi_{\rho}(L,N)|^2}{\lVert \nabla \rho \rVert^2} - \frac{1}{2} \frac{N\Levi_{\rho}(L,L)}{\lVert \nabla \rho \rVert} \le 0$$
		at $p$.
	\end{lem}
	
	\begin{proof}
		By using $L\rho = 0$ and the assumption, we have 
		\begin{align*}
			& \Levi_{\rho^{\eta_2}}(aL + bN, aL + bN)  \\  
			= & \eta_2 \rho^{\eta_2 - 1} \left[ |a|^2 \Levi_{\rho}(L,L) + 2 \text{Re} \left( a\overline{b} \Levi_{\rho}(L,N) \right) + |b|^2\left( \Levi_{\rho}(N,N) + \frac{\eta_2 - 1}{\rho}|N\rho|^2 \right) \right] > 0  
		\end{align*}
		for all $(a,b) \in \CC^2 \setminus (0,0)$ on $\overline{\O}^{\complement} \cap U_p$, and it is equivalent to 
		\begin{equation} \label{1.1}
		|a|^2 \Levi_{\rho}(L,L) - 2|a||b||\Levi_{\rho}(L,N)| + |b|^2\left( \Levi_{\rho}(N,N) + \frac{\eta_2 - 1}{\rho}|N\rho|^2 \right) > 0 
		\end{equation} 
		for all $(a,b) \in \CC^2 \setminus (0,0)$ on $\overline{\O}^{\complement} \cap U_p$. This is because $\eta_2 \rho^{\eta_2 -1} > 0$ and by rotating $a$ or $b$, one can make $2 \text{Re} \left( a\overline{b} \Levi_{\rho}(L,N) \right) = - 2|a||b||\Levi_{\rho}(L,N)| $. We may assume that $\Levi_{\rho}(N,N) + \frac{\eta_2 - 1}{\rho} |N\rho|^2  > 0 $ on $\overline{\O}^{\complement} \cap U_p$, because it blows up as point goes to the boundary. By divding (\ref{1.1}) by $|a|^2$ and letting $x = \frac{|b|}{|a|}$, (\ref{1.1}) is equivalent to
		\begin{equation} \label{1.2}
			\Levi_{\rho}(L,L) - 2x|\Levi_{\rho}(L,N)| + x^2\left( \Levi_{\rho}(N,N) + \frac{\eta_2 - 1}{\rho}|N\rho|^2 \right) > 0 
		\end{equation}
		for all $x \ge 0$ on $\overline{\O}^{\complement} \cap U_p$. The axis of symmetry of the quadratic equation above is $$\frac{|\Levi_{\rho}(L,N)|}{\Levi_{\rho}(N,N) + \frac{\eta_2 - 1}{\rho}|N\rho|^2},$$ 
		which is always positive on $\overline{\O}^{\complement} \cap U_p$. Thus, (\ref{1.2}) if and only if the determinent 
		\begin{equation} \label{1.3}
			\left| \Levi_{\rho}(L,N) \right|^2 - \left( \Levi_{\rho}(L,L) \right) \left( \Levi_{\rho}(N,N) + \frac{\eta_2 - 1}{\rho}|N\rho|^2 \right) < 0 
		\end{equation}
		on $\overline{\O}^{\complement} \cap U_p$. Now, taking a limit to $p$ on (\ref{1.3}) along the real normal direction yields 
		\begin{equation} \label{1.4}
			\left| \Levi_{\rho}(L,N) \right|^2 - (\eta_2 - 1)|N\rho|^2 \lim\limits_{z \rightarrow p}\frac{\Levi_{\rho}(L,L)}{\rho} \le 0 
		\end{equation}
		at $p \in \partial \O$. Here, 
		\begin{equation} \label{1.5}
			\lim\limits_{z \rightarrow p}\frac{\Levi_{\rho}(L,L)(z)}{\rho(z)} = \frac{-(N+\overline{N}) \Levi_{\rho}(L,L)(p)}{-(N+\overline{N})\rho(p) } = \frac{N \Levi_{\rho}(L,L)(p)}{N\rho(p)}.
		\end{equation}
		The explanation of second equality of (\ref{1.5}) is following.
		Since $\Levi_{\rho}(L,L) = 0$ at $p \in \partial \O$, $\Levi_{\rho}(L,L)$ attains the local minimum at $p$ on $\partial \O$. Thus, the directional derivative of $\Levi_{\rho}(L,L)$ along the real tangent vector $N - \overline{N}$ at $p$ is zero, and $(N-\overline{N})\Levi_{\rho}(L,L) = 0$ implies $N\Levi_{\rho}(L,L) = \overline{N}\Levi_{\rho}(L,L)$ at $p$. Finally, since $N\rho =  \frac{\lVert \nabla \rho \rVert}{2}$, (\ref{1.4}) is equivalent to
		
		\begin{equation*}
		\frac{1}{\eta_2-1} \frac{|\Levi_{\rho}(L,N)|^2}{\lVert \nabla \rho \rVert^2} - \frac{1}{2} \frac{N\Levi_{\rho}(L,L)}{\lVert \nabla \rho \rVert} \le 0 
		\end{equation*}
		at $p \in \partial \O$.
	\end{proof}

	\begin{lem} \label{S lem2}
		Let $\Omega \subset\subset \CC^n$ be a pseudoconvex domain with smooth boundary, and $\rho$ be a defining function of $\Omega$. Let $\{ U_{\alpha} \}$ be a chart of $U$ and $L$ be a non-vanishing smooth $(1,0)$ tangent vector field in $U_{\alpha}$ such that $L \rho = 0$. We normalize $L$ as $g(L,L)=\frac{1}{2}$. Fix $\eta_2 \in (1,\infty)$. We denote $N_{\rho}$ by $N$. Define 
		$$ \Sigma_L^{\alpha} := \{ p \in \partial \O \cap U_{\alpha} : \Levi_{\rho}(L,L)(p) = 0  \}. $$
 		If
 		\begin{equation} \label{1.6}
 			\frac{1}{\eta_2-1} \frac{|\Levi_{\rho}(L,N)|^2}{\lVert \nabla \rho \rVert^2} - \frac{1}{2} \frac{N\Levi_{\rho}(L,L)}{\lVert \nabla \rho \rVert} < 0 
 		\end{equation}
 		on $\Sigma_L^{\alpha}$, then there exists a neighborhood $V_L^{\alpha}$ of $\Sigma_L^{\alpha}$ in $U_{\alpha}$ such that
 		\begin{equation*}
 			\Levi_{\rho^{\eta_2}}(aL + bN, aL + bN) > 0
 		\end{equation*}
 		for all $(a,b) \in \CC^2 \setminus (0,0)$ on $\overline{\O}^{\complement} \cap V_L^{\alpha}$.
	\end{lem}
	
	\begin{proof}
		First, note that the assumption (\ref{1.6}) is equivalent to 
		\begin{equation*}
			\left| \Levi_{\rho}(L,N) \right|^2 - (\eta_2 - 1)(N\rho) (N \Levi_{\rho}(L,L)) < 0 
		\end{equation*}
		on $\Sigma_L^{\alpha}$. Now define $F : \O^{\complement} \cap U_{\alpha} \rightarrow \RR$ by \\
		\begin{align*}
			F(z) := & \left| \Levi_{\rho}(L,N)(z) \right|^2 - \left( \Levi_{\rho}(L,L)(z) - \Levi_{\rho}(L,L)(p) \right) \Levi_{\rho}(N,N)(z) \\ 
			& - (\eta_2 - 1) \frac{\Levi_{\rho}(L,L)(z) - \Levi_{\rho}(L,L)(p)}{\rho(z)} |N\rho(z)|^2 
		\end{align*}
		for all $z \in \overline{\O}^{\complement} \cap U_{\alpha}$, where $p \in \partial \O \cap U_{\alpha}$ is the closest point to $z$, and
		\begin{equation*}
			F(z) := \left| \Levi_{\rho}(L,N)(z) \right|^2 - (\eta_2 - 1)(N\rho(z)) (N \Levi_{\rho}(L,L)(z)) 
		\end{equation*}
		for all $z \in \partial \O \cap U_{\alpha}$. Since $\lim\limits_{z \rightarrow p} F(z) = F(p)$ for all $p \in \partial \O \cap U_{\alpha}$, $F$ is continuous on $\O^{\complement} \cap U_{\alpha}$. By the assumption (\ref{1.6}), there exists a neighborhood $V_L^{\alpha}$ of $\Sigma_L^{\alpha}$ in $U_{\alpha}$ such that 
		\begin{align*}
			& \left| \Levi_{\rho}(L,N)(z) \right|^2 - \left( \Levi_{\rho}(L,L)(z) - \Levi_{\rho}(L,L)(p) \right) \Levi_{\rho}(N,N)(z) \\ 
			& - (\eta_2 - 1) \frac{\Levi_{\rho}(L,L)(z) - \Levi_{\rho}(L,L)(p)}{\rho(z)} |N\rho(z)|^2  \\
			= & \left| \Levi_{\rho}(L,N)(z) \right|^2 
			- \left( \Levi_{\rho}(L,L)(z) - \Levi_{\rho}(L,L)(p) \right) \left( \Levi_{\rho}(N,N)(z) + \frac{\eta_2 - 1}{\rho(z)} |N\rho(z)|^2 \right) < 0
		\end{align*}
		for all $z \in \overline{\O}^{\complement} \cap V_L^{\alpha}$. We may assume that $\Levi_{\rho}(N,N)(z) + \frac{\eta_2 - 1}{\rho(z)} |N\rho(z)|^2  > 0 $ for all $z \in \overline{\O}^{\complement} \cap V_L^{\alpha}$, because it blows up as $z$ goes to the boundary. Therefore, the following quadratic function 
		\begin{align*}
			&\left( \Levi_{\rho}(L,L)(z) - \Levi_{\rho}(L,L)(p) \right) -2x \left| \Levi_{\rho}(L,N)(z) \right| \\
			& + x^2 \left( \Levi_{\rho}(N,N)(z) + \frac{\eta_2 - 1}{\rho(z)} |N\rho(z)|^2 \right) > 0
		\end{align*}
		for all $z \in \overline{\O}^{\complement} \cap V_L^{\alpha}$, $x \ge 0$. By letting $x = \frac{|b|}{|a|}$ and multiplying $|a|^2$, we have 
		\begin{align}
		&|a|^2 \left( \Levi_{\rho}(L,L)(z) - \Levi_{\rho}(L,L)(p) \right) -2 |a||b| \left| \Levi_{\rho}(L,N)(z) \right| \label{1.7} \\ 
		&+ |b|^2 \left( \Levi_{\rho}(N,N)(z) + \frac{\eta_2 - 1}{\rho(z)} |N\rho(z)|^2 \right) > 0 \nonumber
		\end{align} 
		for all $z \in \overline{\O}^{\complement} \cap V_L^{\alpha}$, $(a,b) \in \CC^2 \setminus (0,0)$. If $|a| = 0$, then (\ref{1.7}) is equivalent to $\Levi_{\rho}(N,N)(z) + \frac{\eta_2 - 1}{\rho(z)} |N\rho(z)|^2  > 0$, which is automatically satisfied. Now, (\ref{1.7}) implies 
		\begin{align*}
		&|a|^2 \left( \Levi_{\rho}(L,L)(z) - \Levi_{\rho}(L,L)(p) \right) + 2 \Re \left( a\overline{b} \left( \Levi_{\rho}(L,N)(z) \right) \right) \\
		&+ |b|^2 \left( \Levi_{\rho}(N,N)(z) + \frac{\eta_2 - 1}{\rho(z)} |N\rho(z)|^2 \right) > 0
		\end{align*}
		for all $z \in \overline{\O}^{\complement} \cap V_L^{\alpha}$, $(a,b) \in \CC^2 \setminus (0,0)$ because $ \Re \left( a\overline{b} \left( \Levi_{\rho}(L,N)\right) \right) > - |a||b| \left| \Levi_{\rho}(L,N) \right| $. Since $\O$ is pseudoconvex, $\Levi_{\rho}(L,L)(p) \ge 0$ for all $p \in \partial \O$, and thus
		\begin{align}
		&|a|^2 \left( \Levi_{\rho}(L,L)(z) \right) + 2 \Re \left( a\overline{b} \left( \Levi_{\rho}(L,N)(z) \right) \right) \label{ 1.8} \\
		&+ |b|^2 \left( \Levi_{\rho}(N,N)(z) + \frac{\eta_2 - 1}{\rho(z)} |N\rho(z)|^2 \right) > 0 \nonumber
		\end{align}
		for all $z \in \overline{\O}^{\complement} \cap V_L^{\alpha}$, $(a,b) \in \CC^2 \setminus (0,0)$. Finally, by using $L\rho = 0$, (\ref{ 1.8}) is equivalent to 
		\begin{align*}
			\Levi_{\rho^{\eta_2}}(aL + bN, aL + bN)(z) > 0
		\end{align*}
		for all $z \in \overline{\O}^{\complement} \cap V_L^{\alpha}$, $(a,b) \in \CC^2 \setminus (0,0)$.	
	\end{proof}

	\begin{lem} \label{S lem3}
		Let $\Omega \subset\subset \CC^n$ be a pseudoconvex domain with smooth boundary, and $\rho$ be a defining function of $\Omega$. Let $\{ U_{\alpha} \}$ be a chart of $U$ and $L$ be a non-vanishing smooth $(1,0)$ tangent vector field in $U_{\alpha}$ such that $L \rho = 0$. We normalize $L$ as $g(L,L)=\frac{1}{2}$. Fix $\eta_2 \in (1,\infty)$. We denote $N_{\rho}$ by $N$. Define 
		$$ \Sigma_L^{\alpha} := \{ p \in \partial \O \cap U_{\alpha} : \Levi_{\rho}(L,L)(p) = 0  \}. $$
		If
		\begin{equation*}
		\frac{1}{\eta_2-1} \frac{|\Levi_{\rho}(L,N)|^2}{\lVert \nabla \rho \rVert^2} - \frac{1}{2} \frac{N\Levi_{\rho}(L,L)}{\lVert \nabla \rho \rVert} \le 0 
		\end{equation*}
		on $\Sigma_L^{\alpha}$, then there exists a defining function $\widetilde{\rho}$ such that for any small $\nu >0$ so that $\frac{1}{\eta_2 - 1} - \nu > 0$,
		\begin{equation*}
		\frac{1}{\widetilde{\eta}_2-1} \frac{|\Levi_{\widetilde{\rho}}(\widetilde{L},\widetilde{N})|^2}{\lVert \nabla \widetilde{\rho} \rVert^2} - \frac{1}{2} \frac{\widetilde{N}\Levi_{\widetilde{\rho}}(\widetilde{L},\widetilde{L})}{\lVert \nabla \widetilde{\rho} \rVert} < 0 
		\end{equation*}
		on $\Sigma_L^{\alpha}$, where  $\widetilde{\eta}_2 = 1 + \frac{1}{\frac{1}{\eta_2 - 1} - \nu}$ and $\widetilde{N} = N_{\widetilde{\rho}}$. $\widetilde{L}$ is a non-vanishing smooth $(1,0)$ tangent vector field in $U_{\alpha}$ such that $\widetilde{L} = L$ on $\partial \O \cap U_{\alpha}$, and  $\widetilde{L} \widetilde{\rho} = 0$ on $U_{\alpha}$.
	\end{lem}
	
	\begin{proof}
		First, notice that $\frac{1}{\widetilde{\eta}_2 - 1} = \frac{1}{\eta_2 - 1} - \nu$.
		Define $\widetilde{\rho} = \rho e^{\epsilon \psi}$, where $\psi = \norm{z}^2$, $\epsilon$ is a small positive number and we will decide $\epsilon$ later.
		Then by Lemma \ref{keylem}, 
		\begin{align*}
			& \frac{1}{\widetilde{\eta}_2-1} \frac{|\Levi_{\widetilde{\rho}}(\widetilde{L},\widetilde{N})|^2}{\lVert \nabla \widetilde{\rho} \rVert^2} - \frac{1}{2} \frac{\widetilde{N}\Levi_{\widetilde{\rho}}(\widetilde{L},\widetilde{L})}{\lVert \nabla \widetilde{\rho} \rVert}   \\		
			= & \left( \frac{1}{\eta_2-1} - \nu \right) \frac{|\Levi_{\rho}(L,N)|^2}{\lVert \nabla \rho \rVert^2} - \frac{1}{2} \frac{N\Levi_{\rho}(L,L)}{\lVert \nabla \rho \rVert} \\
			& + \left( \frac{1}{\eta_2-1} + 1 - \nu \right) \left( \frac{\epsilon^2}{4}|L\psi|^2 + \epsilon \Re \left[ \frac{\Levi_\rho(L,N)}{\lVert \nabla \rho \rVert} (\overline{L}\psi) \right] \right) - \frac{\epsilon}{8} \\
			\le & \left( \frac{1}{\eta_2-1} - \nu \right) \frac{|\Levi_{\rho}(L,N)|^2}{\lVert \nabla \rho \rVert^2} - \frac{1}{2} \frac{N\Levi_{\rho}(L,L)}{\lVert \nabla \rho \rVert} \\
			& + \left( \frac{1}{\eta_2-1} + 1 - \nu \right) \frac{\epsilon^2}{4}|L\psi|^2 + \frac{\nu}{2} \frac{|\Levi_\rho(L,N)|^2}{\lVert \nabla \rho \rVert^2} + \frac{\left( \frac{1}{\eta_2-1} + 1 - \nu \right)^2}{2\nu} |L\psi|^2 \epsilon^2 - \frac{\epsilon}{8} \\
			\le & \frac{1}{\eta_2-1} \frac{|\Levi_{\rho}(L,N)|^2}{\lVert \nabla \rho \rVert^2} - \frac{1}{2} \frac{N\Levi_{\rho}(L,L)}{\lVert \nabla \rho \rVert} \\
			& - \frac{\nu}{2} \frac{|\Levi_{\rho}(L,N)|^2}{\lVert \nabla \rho \rVert^2} + \left( \frac{1}{\eta_2-1} + 1 - \nu \right) \left( \frac{1}{4} + \frac{\frac{1}{\eta_2-1} + 1 - \nu}{2\nu} \right)|L\psi|^2 \epsilon^2 - \frac{\epsilon}{8}
		\end{align*}
		on $\Sigma_L^{\alpha}$. Now $\frac{1}{\eta_2-1} \frac{|\Levi_{\rho}(L,N)|^2}{\lVert \nabla \rho \rVert^2} - \frac{1}{2} \frac{N\Levi_{\rho}(L,L)}{\lVert \nabla \rho \rVert} \le 0$ by the assumption, and $- \frac{\nu}{2} \frac{|\Levi_{\rho}(L,N)|^2}{\lVert \nabla \rho \rVert^2} \le 0$. We choose sufficiently small $\epsilon > 0$  so that $\left( \frac{1}{\eta_2-1} + 1 - \nu \right) \left( \frac{1}{4} + \frac{\frac{1}{\eta_2-1} + 1 - \nu}{2\nu} \right)|L\psi|^2 \epsilon^2 - \frac{\epsilon}{8} < 0$. All together, we conclude that
	
		\begin{equation*}
			\frac{1}{\widetilde{\eta}_2-1} \frac{|\Levi_{\widetilde{\rho}}(\widetilde{L},\widetilde{N})|^2}{\lVert \nabla \widetilde{\rho} \rVert^2} - \frac{1}{2} \frac{\widetilde{N}\Levi_{\widetilde{\rho}}(\widetilde{L},\widetilde{L})}{\lVert \nabla \widetilde{\rho} \rVert} < 0
		\end{equation*}
		on $\Sigma_L^{\alpha}$.
	\end{proof}

	\begin{lem} \label{S lem4}
		Let $\Omega \subset\subset \CC^n$ be a pseudoconvex domain with smooth boundary, and $\Pi$ be the set of all strongly pseudoconvex points in $\partial \O$. Then for any defining function $\rho$ of $\Omega$ and $\eta_2 \in (1,\infty)$, there exists a neighborhood $V$ of $\Pi$ in $\CC^n$ such that $\rho^{\eta_2}$ is strictly plurisubharmonic on $\overline{\O}^{\complement} \cap V$.
	\end{lem}

	\begin{proof}
		Let $V$ be a neighborhood of $\Pi$ in $\CC^n$, and $\{ V_{\alpha} \}$ be a chart of $V$. We denote $N = N_{\rho}$. It is enough to show that
		\begin{equation} \label{2.155}
			\Levi_{\rho^{\eta_2}}(aL+bN, aL+bN) > 0
		\end{equation}
		for all $(a,b) \in \CC^2 \setminus (0,0)$ on $\overline{\O}^{\complement} \cap V_{\alpha}$, where $L$ is an arbitrary non-vanishing smooth $(1,0)$ tangent vector field in $V_{\alpha}$. As in the proof of Lemma \ref{S lem1}, (\ref{2.155}) is equivalent to 
		\begin{equation} \label{2.156}
			\left| \Levi_{\rho}(L,N) \right|^2 - \left( \Levi_{\rho}(L,L) \right) \left( \Levi_{\rho}(N,N) + \frac{\eta_2 - 1}{\rho}|N\rho|^2 \right) < 0 
		\end{equation}
		on $\overline{\O}^{\complement} \cap V_{\alpha}$. After possibly shrinking $V$, (\ref{2.156}) holds because $\Levi_{\rho}(L,L)>0$ on $\Pi$ and $\Levi_{\rho}(N,N) + \frac{\eta_2 - 1}{\rho}|N\rho|^2$ blows up as point approaches to $\partial \O$. 
	\end{proof}

	\vspace{5mm}
	\begin{proof}[{\bf Proof of Theorem \ref{S main thm}}]
		First, we prove $S(\O) \ge \inf \eta_{\rho}$. For a defining function $\rho$ and $\eta_2 \in (1,\infty)$, if $\rho^{\eta_2}$ is strictly plurisubharmonic on $\overline{\O}^{\complement} \cap U$, then 
		$$ \Levi_{\rho^{\eta_2}}(aL + bN, aL + bN) > 0$$
		for all $(a,b) \in \CC^2 \setminus (0,0)$ on $\overline{\O}^{\complement} \cap U$. By Lemma \ref{S lem1}, this implies that
		$$ \frac{1}{\eta_2-1} \frac{|\Levi_{\rho}(L,N)|^2}{\lVert \nabla \rho \rVert^2} - \frac{1}{2} \frac{N\Levi_{\rho}(L,L)}{\lVert \nabla \rho \rVert} \le 0 $$ 
		on $\Sigma_L$ for all $L$. Therefore, $S(\O) \ge \inf \eta_{\rho}$.
		
		Next, we prove $S(\O) \le \inf \eta_{\rho}$. By Lemma \ref{S lem4}, we only need to consider it in a neighborhood of $\Sigma$ in $\CC^n$. Let $\eta_0 = \inf \eta_{\rho}$. Fix $\widetilde{\eta}_2 > \eta_0$ and choose small $\nu > 0$ so that 
		\begin{align*}
			\eta_2 := 1 + \frac{1}{\frac{1}{\widetilde{\eta}_2 - 1} + \nu}  \phantom{aa} 
			\left( \Leftrightarrow 
			\widetilde{\eta}_2 = 1 + \frac{1}{\frac{1}{\eta_2 - 1} - \nu} \right)
		\end{align*} 
		satisfies $\widetilde{\eta}_2 > \eta_2 > \eta_0$. 
		Let $\{ U_{\alpha} \}$ be a chart of $U$ and $L$ be a non-vanishing smooth $(1,0)$ tangent vector field in $\partial \O \cap U_{\alpha}$ such that $g(L,L)=\frac{1}{2}$. Since $\eta_2 > \eta_0$, there exists $\rho$ such that
		$$ \frac{1}{\eta_2-1} \frac{|\Levi_{\rho}(L,N)|^2}{\lVert \nabla \rho \rVert^2} - \frac{1}{2} \frac{N\Levi_{\rho}(L,L)}{\lVert \nabla \rho \rVert} \le 0 $$  
		on $\Sigma_L^{\alpha}$ for all $L$. Here, we extend $L$ so that $L\rho = 0$ on $U_{\alpha}$. By Lemma \ref{S lem3}, there exists a defining function $\widetilde{\rho}$ such that 
		\begin{equation*}
			\frac{1}{\widetilde{\eta}_2-1} \frac{|\Levi_{\widetilde{\rho}}(\widetilde{L},\widetilde{N})|^2}{\lVert \nabla \widetilde{\rho} \rVert^2} - \frac{1}{2} \frac{\widetilde{N}\Levi_{\widetilde{\rho}}(\widetilde{L},\widetilde{L})}{\lVert \nabla \widetilde{\rho} \rVert} < 0 
		\end{equation*}
		on $\Sigma_L^{\alpha}$ for all $L$. 
		Since $\widetilde{\eta}_2$ is arbitrary, by Lemma \ref{S lem2}, we have $S(\O) \le \inf \eta_{\rho}$.
	\end{proof}

	Together with Lemma \ref{keylem}, we have

	\begin{cor} \label{cor2}
		Let $\Omega \subset\subset \CC^n$ be a pseudoconvex domain with smooth boundary, and $\rho$ be a defining function of $\Omega$. Let $L$ be an arbitrary $(1,0)$ tangent vector field on $\partial \O$. Define
		$$ \Sigma_L := \{ p \in \partial \O : \Levi_{\rho}(L,L)(p) = 0 \}. $$
		Let $\psi$ be a smooth function defined on $U$, and $\eta_{\psi} $ be the infimum of $\eta_2 \in (1,\infty)$ satisfying 
		$$ \left( \frac{1}{\eta_2-1} + 1 \right) \left| \frac{\Levi_{\rho}(L,N)}{\lVert \nabla \rho \rVert} + \frac{1}{2}(L\psi) \right|^2 -
		\left( \frac{|\Levi_{\rho}(L,N)|^2}{\lVert \nabla \rho \rVert^2} + \frac{1}{2} \frac{N\Levi_{\rho}(L,L)}{\lVert \nabla \rho \rVert} \right) - \frac{1}{4}\Levi_{\psi}(L,L)  \le 0 $$
		on $\Sigma_L$ for all $L$. Here, $N = N_{\rho}$ and we extend $L$ so that $L \rho = 0$ on $U$.
		Then 
		$$S(\O) = \inf \text{ } \eta_{\psi},$$
		where the infimum is taken over all smooth functions $\psi$.
	\end{cor}

	\vspace{5mm}
	

	\vspace{5mm}
	\section{\bf Strong Stein neighborhood basis} \label{section strong Stein neighborhood basis}
	
	In 2012, Sahuto\u{g}lu (\cite{Sahutoglu2012}) gave several characterizations for $\overline{\O}$ to have a strong Stein neighborhood basis. In this section, using one of the characterizations, we prove that the existence of the Steinness index and the existence of a strong Stein neighborhood basis are equivalent.

	\begin{thm} \label{Sindex eqi SSnbhd}
		Let $\O \subset \subset \CC^n$ be a pseudoconvex domain with smooth boundary. Then $S(\O)$ exists if and only if $\overline{\O}$ has a strong Stein neighborhood basis.
	\end{thm}
	
	\begin{proof}
		First, assume $S(\O)$ exists. Then there exists a defining function $\rho$ of $\O$ and $\eta_2 \in (1, \infty)$ such that $\rho^{\eta_2}$ is strictly plurisubharmonic on $\overline{\O}^{\complement} \cap U$. Since level sets of $\rho^{\eta_2}$ and $\rho$ are same, $\rho$ is the desired defining function. 
		
		Now, suppose that $\overline{\O}$ has a strong Stein neighborhood basis. Then by Sahuto\u{g}lu (\cite{Sahutoglu2012}), there exist a defining function $\rho$ of $\O$ and $c > 0$ such that 
		\begin{equation*}
		\Levi_{\rho}(L,L) \ge c \rho \lVert L \rVert^2
		\end{equation*}
		on $\overline{\O}^{\complement} \cap U$, where $L$ is a smooth $(1,0)$ tangent vector field in $U$ with $L \rho = 0$. By letting $\norm{L}^2 = \frac{1}{2}$ and dividing it by $\rho$, we have
		\begin{equation*} 
		\frac{\Levi_{\rho}(L,L)}{\rho} \ge \frac{c}{2}
		\end{equation*}
		Define $\Sigma_L := \{ p \in \Sigma : \Levi_{\rho}(L,L)(p) = 0 \}$ and denote $N = N_{\rho}$.  Taking a limit to $p \in \Sigma_L$ along the real normal direction gives
		\begin{align*}
		\frac{1}{2} \frac{N \Levi_{\rho}(L,L)}{\lVert \nabla \rho \rVert} \ge \frac{c}{8}
		\end{align*}
		on $\Sigma_L$. On the other hands, since $\Sigma$ and $ \{ L \in T^{1,0}_p(\partial \O) : \norm{L}^2 = \frac{1}{2} \} $ are compact,
		$$ \frac{|\Levi_{\rho}(L,N)|^2}{\lVert \nabla \rho \rVert^2} $$
		has the maximum value $M \ge 0$ on $\Sigma$. Therefore, for $\eta_2 \in (1,\infty)$
		\begin{align*}
		\frac{1}{\eta_2-1} \frac{|\Levi_{\rho}(L,N)|^2}{\lVert \nabla \rho \rVert^2} - \frac{1}{2} \frac{N\Levi_{\rho}(L,L)}{\lVert \nabla \rho \rVert} \le 
		\frac{M}{\eta_2 - 1} - \frac{c}{8}
		\end{align*}
		on $\Sigma_L$. By choosing $\eta_2 > \frac{8M}{c} + 1$, we have
		\begin{align*}
		\frac{1}{\eta_2-1} \frac{|\Levi_{\rho}(L,N)|^2}{\lVert \nabla \rho \rVert^2} - \frac{1}{2} \frac{N\Levi_{\rho}(L,L)}{\lVert \nabla \rho \rVert} \le 0
		\end{align*}
		on $\Sigma_L$. By Theorem \ref{S main thm}, we proved $S(\O)$ exists.
	\end{proof}

	\vspace{5mm}
	
	\section{\bf Steinness index of worm domains} \label{section Steinness index of worm domains}

	In this section, we calculate the exact value of the Steinness index of worm domains, and prove Theorem \ref{main}. Recall first the definition of worm domains. 
	
	\begin{defn} \label{worm defn}
		The {\it worm domain} $\Omega_{\beta}$ $(\beta > \frac{\pi}{2})$ is defined by
		$$ \Omega_{\beta} := \left\lbrace (z,w)\in \CC^2 : \rho(z,w) = \left| z - e^{i \log|w|^2} \right|^2 - (1 - \phi(\log|w|^2) ) < 0  \right\rbrace $$
		where $\phi : \RR \rightarrow \RR$ is a fixed smooth function with the following properties : \\
		1. $\phi(x) \ge 0$, $\phi$ is even and convex. \\
		2. $\phi^{-1}(0) = I_{\beta - \frac{\pi}{2}} = [-(\beta - \frac{\pi}{2}), \beta - \frac{\pi}{2} ].$ \\
		3. $\exists$ $a>0$ such that $\phi(x)>1$ if $x<-a$ or $x>a$.\\
		4. $\phi'(x) \neq 0$ if $\phi(x) = 1$.
	\end{defn}

	Let $\rho(z,w) = \left| z - e^{i \log|w|^2} \right|^2 - (1 - \phi(\log|w|^2) )$, and $U$ be a neighborhood of $\partial \O_{\beta}$.
	Define
	\begin{align*}
	& L = \frac{1}{\sqrt{\left| \frac{\partial \rho}{\partial z} \right|^2 + \left| \frac{\partial \rho}{\partial w} \right|^2 }} \left( \frac{\partial \rho}{\partial w} \frac{\partial}{\partial z} - \frac{\partial \rho}{\partial z} \frac{\partial}{\partial w} \right), \\ 
	& N = \frac{1}{\sqrt{\left| \frac{\partial \rho}{\partial z} \right|^2 + \left| \frac{\partial \rho}{\partial w} \right|^2 }} \left( \frac{\partial \rho}{\partial \overline{z}} \frac{\partial}{\partial z} + \frac{\partial \rho}{\partial \overline{w}} \frac{\partial}{\partial w} \right).
	\end{align*}
	Let $\Sigma$ be the set of all weakly pseudoconvex points in $\partial \O_{\beta}$. Then 
	\begin{equation*}
	\Sigma = \left\lbrace (0,w) \in \CC^2 : \left| \log|w|^2 \right| \le \beta - \frac{\pi}{2}  \right\rbrace.
	\end{equation*}
	By direct calculation, we have
	\begin{align*}
	\lVert \nabla \rho \rVert = 2, \phantom{aaa}  
	L = e^{-i\log|w|^2} \frac{\partial}{\partial w}, 	\phantom{aaa}
	N = - e^{i\log|w|^2} \frac{\partial}{\partial z},
	\end{align*}
	\begin{align*}
	\Levi_{\rho}(L,N) = \frac{i}{w}e^{-i\log|w|^2}, \phantom{aaa}
	N\Levi_{\rho}(L,L) = -\frac{1}{|w|^2}
	\end{align*}
	on $\Sigma$.

	\begin{lem} [Riccati equations \cite{Liu2017}] \label{Riccati}
		For $a,b >0$ and $t>0$, the following Riccati euqation 
		$$ \frac{d}{dt}s(t) = a (s(t))^2 - \frac{s(t)}{t} + \frac{b}{t^2} $$
		has the solution 
		\begin{equation*}
		s(t) = - \sqrt{\frac{b}{a}} \frac{\cot(\sqrt{ab}\log t + \phi )}{t}
		\end{equation*} 
		for arbitrary $\phi$.
	\end{lem}

	\begin{thm} \label{S worm}
		Let $\Omega_{\beta}$ $(\beta > \frac{\pi}{2})$ be a worm domain. Then 
		\begin{equation*}
		S(\O_{\beta}) = 
		\begin{cases}
		\frac{\pi}{2(\pi - \beta)}  & \text{ for } \frac{\pi}{2} < \beta < \pi  \\
		\infty 						& \text{ for } \pi \le \beta
		\end{cases}
		\end{equation*}
	\end{thm}

	\begin{proof}
		We will use Corollary \ref{cor2}. Let $\alpha = \frac{1}{\eta_2-1} + 1$ for $\eta_2 \in (1, \infty)$. Suppose that there exists a smooth function $\psi$ defined on $U$ such that 
		\begin{equation}
		\alpha \left| \frac{\Levi_{\rho}(L,N)}{\lVert \nabla \rho \rVert} + \frac{1}{2}(L\psi) \right|^2 -
		\left( \frac{|\Levi_{\rho}(L,N)|^2}{\lVert \nabla \rho \rVert^2} + \frac{1}{2} \frac{N\Levi_{\rho}(L,L)}{\lVert \nabla \rho \rVert} \right) - \frac{1}{4}\Levi_{\psi}(L,L)  \le 0  \label{3.1}
		\end{equation}
		on $\Sigma$. By the calculation above, 
		\begin{equation*}
		\frac{|\Levi_{\rho}(L,N)|^2}{\lVert \nabla \rho \rVert^2} + \frac{1}{2} \frac{N\Levi_{\rho}(L,L)}{\lVert \nabla \rho \rVert} = \frac{1}{4|w|^2} - \frac{1}{4|w|^2} = 0
		\end{equation*}
		on $\Sigma$. Consequently, (\ref{3.1}) is equivalent to 
		\begin{equation}
		\alpha \left| \frac{i}{w} + \frac{\partial \psi}{\partial w} \right|^2 - \frac{\partial^2 \psi}{\partial w \partial \overline{w}} \le 0 \label{3.2}
		\end{equation}
		on $\Sigma$. Let $w = re^{i\theta}$. Then (\ref{3.2}) implies
		\begin{align}
		& \alpha \left( \frac{\partial \psi}{\partial w} \cdot \frac{\partial \psi}{\partial \overline{w}} + \frac{1}{|w|^2} + 2 \Re \left( \frac{\partial \psi}{\partial \overline{w}} \cdot \frac{i}{w} \right) \right) -  \frac{\partial^2 \psi}{\partial w \partial \overline{w}} \nonumber \\
		= & \frac{\alpha}{4} \psi_r^2 + \frac{\alpha}{4r^2} \psi_{\theta}^2 + \frac{\alpha}{r^2} - \frac{\alpha}{r^2} \psi_{\theta} - \frac{1}{4} \psi_{rr} - \frac{1}{4r} \psi_r - \frac{1}{4r^2} \psi_{\theta \theta} \le 0 \label{3.3}
		\end{align}
		on $\Sigma$. Notice that $\int_{0}^{2\pi} \psi_{\theta} d\theta = 0$, $\int_{0}^{2\pi} \psi_{\theta \theta} d\theta = 0$ and $\int_{0}^{2\pi} \psi_{\theta}^2 d\theta \ge 0$. Also by Schwarz's lemma $ \left( \int_{0}^{2\pi} \psi_r d\theta \right)^2 \le \left( \int_{0}^{2\pi} d\theta \right) \left( \int_{0}^{2\pi} \psi_r^2 d\theta \right)$.
		Thus integrating on the both sides of (\ref{3.3}) with respect to $\theta$ gives
		\begin{align}
		\frac{\alpha}{8\pi} \left( \int_{0}^{2\pi} \psi_r d\theta \right)^2 + \frac{2\pi \alpha}{r^2} - \frac{1}{4} \int_{0}^{2\pi} \psi_{rr} d\theta - \frac{1}{4r} \int_{0}^{2\pi} \psi_r d\theta \le 0 \label{3.4}		
		\end{align}
		for all $r \in [e^{- \left( \frac{\beta}{2} -\frac{\pi}{4} \right) } , e^{\frac{\beta}{2} -\frac{\pi}{4}} ]$.  Define $\Psi(r) := \frac{1}{2\pi} \int_{0}^{2\pi} \psi(r,\theta) d\theta$. Then (\ref{3.4}) is equivalent to 
		\begin{align*}
		\frac{\pi \alpha}{2} \Psi_r^2 + \frac{2\pi \alpha}{r^2} - \frac{\pi}{2} \Psi_{rr} - \frac{\pi}{2r} \Psi_r \le 0 		
		\end{align*}
		for all $r \in [e^{- \left( \frac{\beta}{2} -\frac{\pi}{4} \right) } , e^{\frac{\beta}{2} -\frac{\pi}{4}}]$, where $\Psi_r(r) = \frac{d}{dr}\Psi(r)$. Letting $s(r) = \Psi_r(r)$, we have
		\begin{align*}
		- s' + \alpha s^2 - \frac{s}{r} + \frac{4 \alpha}{r^2} \le 0
		\end{align*}
		for all $r \in [e^{- \left( \frac{\beta}{2} -\frac{\pi}{4} \right) } , e^{\frac{\beta}{2} -\frac{\pi}{4}}]$. Suppose $s(1) = s_0$. By the comparison principle of ordinary differential equation, and Lemma \ref{Riccati} 
		\begin{align*}
		s(r) \ge - 2 \frac{\cot(2\alpha \log r + \phi_0)}{r}
		\end{align*}
		for all $r \in [e^{- \left( \frac{\beta}{2} -\frac{\pi}{4} \right) } , e^{\frac{\beta}{2} -\frac{\pi}{4}}]$, where $\phi_0$ is a constant such that
		\begin{align*}
		s_0 = - 2 \cot(\phi_0).
		\end{align*} 
		Since $s(r)$ is a smooth function on $[e^{- \left( \frac{\beta}{2} -\frac{\pi}{4} \right) } , e^{\frac{\beta}{2} -\frac{\pi}{4}}]$, the period of the cotangent function is $\pi$, and $- \alpha \left( \beta - \frac{\pi}{2} \right) \le 2\alpha \log r \le \alpha \left( \beta - \frac{\pi}{2} \right)$, the following must hold.
		\begin{align*} 
		\alpha \left( \beta - \frac{\pi}{2} \right) < \frac{\pi}{2} 
		\end{align*}
		which is equivalent to 
		\begin{align} \label{must}
		\frac{1}{\eta_2 - 1} < \frac{2(\pi - \beta)}{2\beta - \pi}.
		\end{align}
		If $\beta \ge \pi$, then (\ref{must}) never holds. This proves that $S(\O_{\beta}) = \infty$ for all $\beta \ge \pi$. 
		
		Now we assume $\frac{\pi}{2} < \beta < \pi$. Then (\ref{must}) is equivalent to 
		\begin{align}
		\eta_2 > \frac{\pi}{2(\pi-\beta)}.
		\end{align}
		This shows that there does not exist any smooth function $\psi$ defined on $U$ satisfying (\ref{3.1}) if $\eta_2 \le \frac{\pi}{2(\pi-\beta)}$. 
		Therefore, $S(\O_{\beta}) \ge \frac{\pi}{2(\pi-\beta)}$.
		Next, we prove that there exists a smooth function $\psi$ defined on $U$ satisfying (\ref{3.1}) if $\eta_2 > \frac{\pi}{2(\pi-\beta)}$. It is sufficient to find a smooth function $\psi$ on $\Sigma$ because we can extend $\psi$ to $U$. 
		By the argument above, (\ref{3.1}) is equivalent to 
		\begin{align} \label{3.7}
		\frac{\alpha}{4} \psi_r^2 + \frac{\alpha}{4r^2} \psi_{\theta}^2 + \frac{\alpha}{r^2} - \frac{\alpha}{r^2} \psi_{\theta} - \frac{1}{4} \psi_{rr} - \frac{1}{4r} \psi_r - \frac{1}{4r^2} \psi_{\theta \theta} \le 0.
		\end{align}
		We assume $\psi(r, \theta)$ is independent of $\theta$ and let $s(r) = \psi_r$. Then (\ref{3.7}) becomes
		\begin{align} \label{3.8}
		- s' + \alpha s^2 - \frac{s}{r} + \frac{4 \alpha}{r^2} \le 0.
		\end{align}
		By Lemma \ref{Riccati}, $- s' + \alpha s^2 - \frac{s}{r} + \frac{4 \alpha}{r^2} = 0$ has a solution 
		\begin{align*}
		s(r) = - 2 \frac{\cot(2\alpha \log r + \frac{\pi}{2})}{r}.
		\end{align*}
		Since $\eta_2 > \frac{\pi}{2(\pi-\beta)}$, $s(r)$ is a well-defined smooth function on $[e^{- \left( \frac{\beta}{2} -\frac{\pi}{4} \right) } , e^{\frac{\beta}{2} -\frac{\pi}{4}} ]$ and hence $\psi(r,\theta) = \int s(r) dr$ satisfies (\ref{3.1}) on $\Sigma$. Therefore, by Corollary \ref{cor2}, $S(\O_{\beta}) = \frac{\pi}{2(\pi - \beta)}$.
	\end{proof}

	\vspace{5mm}
	\begin{proof}[{\bf Proof of Theorem \ref{main}}]
		
		Notice that the second and fourth conditions are equivalent by Theorem \ref{Sindex eqi SSnbhd}. We will show that the ranges of $\beta$ for each condition are same.
		First, Liu (\cite{Liu2017}) showed $DF(\O_{\beta}) = \frac{\pi}{2\beta}$, which implies that $\frac{1}{2} < DF(\O_{\beta}) < 1$ if and only if $\frac{\pi}{2} < \beta < \pi$. 
		Next, by Theorem \ref{S worm}, the existence of $S(\O_{\beta})$ is equivalent to $\frac{\pi}{2} < \beta < \pi$.
		Finally, one can prove that the third condition is equivalent to $\frac{\pi}{2} < \beta < \pi$ using Theorem 5.1 and 5.2 in \cite{BedForn1978} (see also section 5 in \cite{BedForn1978}). 
		If one of the conditions holds, then $DF(\O_{\beta}) = \frac{\pi}{2\beta}$ and $S(\O_{\beta}) = \frac{\pi}{2(\pi-\beta)}$ implies the last equality.
	\end{proof}

	\vspace{5mm}

	\section{\bf Steinness index of finite type domains} \label{section Steinness index of finite type domains}

	Theorem \ref{S main thm} says that the Steinness index is characterized by some differential inequality on the set of all weakly pseudoconvex boundary points $\Sigma$. Here, we show that considering the set $\Sigma_{\infty}$ of infinite type boundary points suffices to characterize the Steinness index (Theorem \ref{S main thm infinit type}).

	\begin{lem} [\cite{Sahutoglu2012}]  \label{finite type lem}
		Let $\O \subset \subset \CC^n$ be a domain with smooth boundary, and $K$ be a compact subset of $\partial \O$. Assume that $z$ is of finite type for every $z \in K$ and $h \in C^{\infty}(\overline{\O})$ is given. Then for every $j > 0$ there exists $h_j \in C^{\infty}(\overline{\O})$ such that $|h_j - h| \le \frac{1}{j}$ uniformly on $\overline{\O}$ and $\Levi_{h_j}(X,X) \ge j \norm{X}^2$ for all $X \in T_p^{1,0}(\CC^n)$ on $K$.
	\end{lem}

	\begin{lem} \label{finite type lem2}
		Let $\O \subset \subset \CC^n$ be a domain with smooth boundary, and $\Sigma_{\infty} \subset \partial \O$ be the set of all infinite type boundary points. Assume that there exist a neighborhood $V$ of $\Sigma_{\infty}$ in $\CC^n$, a defining function $\rho$ of $\O$, and $\eta_2 > 1$ such that $\rho^{\eta_2}$ is strictly plurisubharmonic on $\overline{\O}^{\complement} \cap V$. Then there exists a defining function $\widetilde{\rho}$ such that $\widetilde{\rho}^{\eta_2}$ is strictly plurisubharmonic on $\overline{\O}^{\complement} \cap U$.
	\end{lem}
	
	\begin{proof} [{\it Proof} {\rm (Based on \cite{Sahutoglu2012})}]
		Let $\Sigma_f$ be the set of all finite type boundary points, and $\Pi$ be the set of all strongly pseudoconvex boundary points. Let $\Sigma_0 := \Sigma_f \setminus \Pi$.  If $\Sigma_0 = \emptyset$, then by Lemma \ref{S lem4}, $\rho^{\eta_2}$ is strictly plurisubharmonic on $\overline{\O}^{\complement} \cap U$. Assume that $\Sigma_0$ is non-empty.
		
		Let $\chi : \RR \rightarrow \RR$ be a smooth increasing convex function such that $\chi(t) = 0$ if $t \le 0$ and $\chi(t) > 0$ if $t > 0$. Let $\beta = \frac{1}{\eta_2 - 1} + 1$ and $A := \max \{ 4 + \beta \chi'(t) : 0 \le t \le 2 \}$. Then by Lemma \ref{finite type lem}, there exist a sequence of $\phi_j \in C^{\infty}(\overline{\O})$ and neighborhoods $V_1, V_2, V_3$ in $\CC^n$ with $\Sigma_{\infty} \subset\subset V_1 \subset\subset V_2 \subset\subset V_3 \subset\subset V$ such that : \\
		\begin{itemize}
			\item[1.] $ -6 \ln A < \phi_j < - \ln A$ on $\overline{\O}$. 
			\item[2.] $ -6 \ln A < \phi_j < - 3\ln A$ on $V_2 \cap \partial \O$.
			\item[3.] $ -2 \ln A < \phi_j < - \ln A$ on $\partial \O \setminus V_3$.
			\item[4.] $\Levi_{\phi_j}(X,X) > j A^6 \norm{X}^2$ for all $X \in T_p^{1,0}(\CC^n)$ on $\partial \O \setminus V_1$. \\
		\end{itemize}
		Let $h_j = e^{\phi_j}$, $a = A^{-3}$, $\chi_a(t) = \chi(t-a)$, $\psi_j = \chi_a \circ h_j$ and $\widetilde{\rho}_j = \rho e^{\psi_j}$. Then $\psi_j \equiv 0$ on a neighborhood $U_1$ of $\overline{V_2 \cap \partial \O}$ in $\CC^n$, hence $\widetilde{\rho_j}^{\eta_2}$ is strictly plurisubharmonic on $\overline{\O}^{\complement} \cap U_1$. Also,
		\begin{align}
			& \Levi_{h_j}(X,X) 
			= e^{\phi_j} \Levi_{\phi_j}(X,X) + e^{\phi_j}|X\phi_j|^2 \\
			= & e^{\phi_j} \Levi_{\phi_j}(X,X) + e^{-\phi_j}|X h_j|^2 
			> j \norm{X}^2 + A|X h_j|^2 \nonumber
		\end{align}
		for all $X \in T_p^{1,0}(\CC^n)$ on $\partial \O \setminus V_1$.

		Let $\{ U_{\alpha} \}$ be a chart of $U$. Let $\widetilde{L}, L$ be non-vanishing smooth $(1,0)$ tangent vector fields in $U_{\alpha}$ such that $\widetilde{L}=L$ on $\partial \O \cap U_{\alpha}$, and $\widetilde{L} \widetilde{\rho}_j = 0$, $L \rho = 0$ on $U_{\alpha}$. We normalize $\widetilde{L}, L$ as $\lVert \widetilde{L} \rVert^2 = \norm{L}^2 = \frac{1}{2}$. Denote $\widetilde{N} = N_{\widetilde{\rho}}$, $N = N_{\rho}$. 
		Define
		$$ \Sigma_L^{\alpha} := \{ p \in (\Sigma_0 \setminus V_1) \cap U_{\alpha} : \Levi_{\rho}(L,L)(p) = 0 \}. $$
		Let $F(\rho, \eta_2) = \frac{1}{\eta_2 -1} \frac{|\Levi_{\rho}(L,N)|^2}{\lVert \nabla \rho \rVert^2} - \frac{1}{2} \frac{N\Levi_{\rho}(L,L)}{\lVert \nabla \rho \rVert}$.
		Then by using Lemma \ref{keylem}, on $\Sigma_L^{\alpha}$,
		\begin{align} \label{6.2}
			F(\widetilde{\rho}_j, \eta_2) 
			\le F(\rho, \eta_2) + \beta \frac{|\Levi_{\rho}(L,N)|}{\lVert \nabla \rho \rVert}|L\psi_j| + \frac{\beta}{4}|L\psi_j|^2 
			 - \frac{1}{4}\Levi_{\psi_j}(L,L).
		\end{align}
		Here,
		\begin{align*}
			\beta \frac{|\Levi_{\rho}(L,N)|}{\lVert \nabla \rho \rVert}|L\psi_j| 
			= \beta \chi_a'(h_j) \frac{|\Levi_{\rho}(L,N)|}{\lVert \nabla \rho \rVert}|Lh_j|
			\le \chi_a'(h_j) \left( \frac{\beta^2 |\Levi_{\rho}(L,N)|^2}{4\lVert \nabla \rho \rVert^2} + |Lh_j|^2 \right),
		\end{align*}
		\begin{align*}
			\frac{\beta}{4}|L\psi_j|^2 = \frac{\beta}{4}(\chi_a'(h_j))^2 |Lh_j|^2,
		\end{align*}
		\begin{align*}
			- \frac{1}{4}\Levi_{\psi_j}(L,L) 
			&= -\frac{1}{4} \chi_a'(h_j) \Levi_{h_j}(L,L) - \frac{1}{4} \chi_a''(h_j) |Lh_j|^2 \\
			&< -\frac{\chi_a'(h_j)}{8}j - \frac{\chi_a'(h_j)}{4}A|Lh_j|^2 - \frac{\chi_a''(h_j)}{4}|Lh_j|^2.
		\end{align*}
		All together, the right-hand side of (\ref{6.2}) is negative if and only if
		\begin{equation} \label{6.3}
			\frac{F(\rho, \eta_2)}{\chi_a'(h_j)} + \frac{\beta^2 |\Levi_{\rho}(L,N)|^2}{4\lVert \nabla \rho \rVert^2} 
			+ \left( 1 + \frac{\beta \chi_a'(h_j)}{4} - \frac{A}{4} - \frac{\chi_a''(h_j)}{4\chi_a'(h_j)} \right)|Lh_j|^2 - \frac{j}{8} < 0.
		\end{equation}
		Since $h_j - a < A^{-1} - A^{-3} < 2$ on $\overline{\O}$, $4 + \beta \chi_a'(h_j) \le A$ by the definition of $A$. Since $\partial \O$ is compact, $\frac{\beta^2 |\Levi_{\rho}(L,N)|^2}{4\lVert \nabla \rho \rVert^2} $ is bounded. Since $\rho^{\eta_2}$ is strictly plurisubharmonic on $V$, $F(\rho, \eta_2) \le 0$ on $\partial \O \cap V$ by Lemma \ref{S lem1}. For the outside of $V$, $F(\rho, \eta_2)$ may be positive, but since $\chi_a'(h_j) > c$ on $\partial \O \setminus V_3$ for some positive constant $c>0$, $\frac{F(\rho, \eta_2)}{\chi_a'(h_j)}$ is bounded. Therefore, there exists a sufficiently large $j>0$ such that (\ref{6.3}) holds. Hence, $F(\widetilde{\rho}_j, \eta_2) < 0$ on $\Sigma_L^{\alpha}$ for all $\alpha$ and $\widetilde{L}$. By Lemma \ref{S lem2}, there exists a neighborhood $U_2$ of $\Sigma_0 \setminus V_1$ such that $\widetilde{\rho_j}^{\eta_2}$ is strictly plurisubharmonic on $\overline{\O}^{\complement} \cap U_2$. Finally, by Lemma \ref{S lem4}, there exists a neighborhood $U_3$ of $\Pi$ such that $\widetilde{\rho_j}^{\eta_2}$ is strictly plurisubharmonic on $\overline{\O}^{\complement} \cap U_3$. Since $U_1 \cup U_2 \cup U_3$ is a neighborhood of $\partial \O$, the proof is completed.
	\end{proof}
	
	\begin{thm} \label{S main thm infinit type} 
		Let $\Omega \subset\subset \CC^n$ be a pseudoconvex domain with smooth boundary, and $\rho$ be a defining function of $\Omega$. Let $L$ be an arbitrary $(1,0)$ tangent vector field on $\partial \O$, and $\Sigma_{\infty}$ be the set of all infinite type boundary points. Define
		$$ \Sigma_{\infty,L} := \{ p \in \Sigma_{\infty} : \Levi_{\rho}(L,L)(p) = 0 \}. $$
		Let $\eta_{\rho} $ be the infimum of $\eta_2 \in (1,\infty)$ satisfying 
		$$ \frac{1}{\eta_2-1} \frac{|\Levi_{\rho}(L,N)|^2}{\lVert \nabla \rho \rVert^2} - \frac{1}{2} \frac{N\Levi_{\rho}(L,L)}{\lVert \nabla \rho \rVert} \le 0 $$
		on $\Sigma_{\infty, L}$ for all $L$. Here, $N = N_{\rho}$ and we extend $L$ so that $L \rho = 0$ on $U$.
		Then 
		$$S(\O) = \inf \text{ } \eta_{\rho}$$
		where the infimum is taken over all smooth defining functions $\rho$.
	\end{thm}
	\begin{proof}
		The proof is same as that of Theorem \ref{S main thm} except that we use Lemma \ref{finite type lem2} instead of Lemma \ref{S lem4}.
	\end{proof}
	
	\begin{cor} \label{finite type}
		Let $\Omega \subset\subset \CC^n$ be a pseudoconvex domain with smooth boundary. Assume that all boundary points of $\O$ are of finite type. Then $S(\O)=1$.
	\end{cor}

	\begin{cor} \label{finite infinite type points}
		Let $\Omega \subset\subset \CC^n$ be a pseudoconvex domain with smooth boundary. Assume that the set of all infinite type boundary points of $\O$ is finite. Then $S(\O)=1$.
	\end{cor}
	\begin{proof}
		Suppose that the set of all infinite type boundary points is one point, says $\Sigma_{\infty} = \{ p_0 \}$.
		We denote $\widetilde{N} = N_{\widetilde{\rho}}$, $N = N_{\rho}$. Let $\widetilde{L}$, $L$ be $(1,0)$ tangent vector fields on $U$ such that $\widetilde{L} = L$ on $\partial \O$, $\widetilde{L} \widetilde{\rho} = 0$, $L\rho = 0$ on $U$, and $\Levi_{\rho}(L,L) = 0$ at $p_0$.
		First, by Lemma \ref{keylem}, 
		\begin{equation} \label{6.4}
			\frac{|\Levi_{\widetilde{\rho}}(\widetilde{L}, \widetilde{N})|^2}{\lVert \nabla \widetilde{\rho} \rVert^2} + \frac{1}{2} \frac{\widetilde{N} \Levi_{\widetilde{\rho}}(\widetilde{L},\widetilde{L})}{\lVert \nabla \widetilde{\rho} \rVert} 
			= \frac{|\Levi_{\rho}(L,N)|^2}{\lVert \nabla \rho \rVert^2} + \frac{1}{2} \frac{N\Levi_{\rho}(L,L)}{\lVert \nabla \rho \rVert} + \frac{1}{4}\Levi_{\phi}(L,L),
		\end{equation}
		where $\widetilde{\rho} = \rho e^{\phi}$.  By letting $\phi = \alpha \norm{z}^2$ and choosing $\alpha>0$ sufficiently large, we can make the left-hand side of (\ref{6.4}) positive at $p_0$ for all $\widetilde{L}$. Hence, we may assume that there exists a defining function $\rho$ of $\O$ such that
		\begin{equation} \label{6.5}
			\frac{|\Levi_{\rho}(L,N)|^2}{\lVert \nabla \rho \rVert^2} + \frac{1}{2} \frac{N\Levi_{\rho}(L,L)}{\lVert \nabla \rho \rVert} > 0
		\end{equation}
		for all $L$ at $p_0$.
		Now, if $L(z) = \sum_{j=1}^{n} a_j(z) \frac{\partial}{\partial z_j}$, then $-2\frac{\Levi_{\rho}(L,N)}{\norm{\nabla \rho}}$ can be represented by
		$$ -2\frac{\Levi_{\rho}(L,N)}{\norm{\nabla \rho}}(z) = \sum_{j=1}^{n} b_j(z) a_j(z) $$
		for some complex-valued functions $b_j(z)$. Define $\psi(z) = \sum_{j=1}^{n} b_j(p_0) z_j$. Then
		$$ L \psi(p_0) = -2\frac{\Levi_{\rho}(L,N)}{\norm{\nabla \rho}}(p_0), \phantom{aaaa} \Levi_{\psi}(L,L)(p_0) = 0.$$ 
		If $\widetilde{\rho} = \rho e^{\psi}$, then by Lemma \ref{keylem} and (\ref{6.5}) at $p_0$
		\begin{align} \label{6.6}
			&\frac{1}{\eta_2 - 1}\frac{|\Levi_{\widetilde{\rho}}(\widetilde{L}, \widetilde{N})|^2}{\lVert \nabla \widetilde{\rho} \rVert^2} 
			+ \frac{1}{2} \frac{\widetilde{N} \Levi_{\widetilde{\rho}}(\widetilde{L},\widetilde{L})}{\lVert \nabla \widetilde{\rho} \rVert} \\ \nonumber
			= &\left( \frac{1}{\eta_2-1} + 1 \right) \left| \frac{\Levi_{\rho}(L,N)}{\lVert \nabla \rho \rVert} + \frac{1}{2}(L\psi) \right|^2 -
			\left( \frac{|\Levi_{\rho}(L,N)|^2}{\lVert \nabla \rho \rVert^2} + \frac{1}{2} \frac{N\Levi_{\rho}(L,L)}{\lVert \nabla \rho \rVert} \right) - \frac{1}{4}\Levi_{\psi}(L,L) \\ \nonumber
			= & -\left( \frac{|\Levi_{\rho}(L,N)|^2}{\lVert \nabla \rho \rVert^2} + \frac{1}{2} \frac{N\Levi_{\rho}(L,L)}{\lVert \nabla \rho \rVert} \right) < 0 \nonumber
		\end{align}
		for all $L$ and $\eta_2 > 1$. Therefore, by Theorem \ref{S main thm infinit type}, $S(\O)=1$. If the number of points in $\Sigma_{\infty}$ is more than 1, then one can construct a smooth function $\psi$ on $U$ satisfying (\ref{6.6}) at all infinite type points using similar argument as above.
	\end{proof}

	\vspace{5mm}
	
	\section{\bf Steinness index of convex domains} \label{section Steinness index of convex domains}

	Assume that $\O \subset \subset \CC^n$ is a domain with $C^k (k \ge 1)$-smooth boundary. The primary goal is to prove that if $\O$ is convex then $DF(\O) = 1$ and $S(\O)=1$. In fact, we prove more: there exists a defining function which is strictly convex on $\CC^n \setminus \partial \O$. For this, we need the following notion of {\it smooth maximum}. For $\epsilon > 0$, let $\chi_{\epsilon} : \RR \rightarrow \RR$ be a smooth function such that $\chi_{\epsilon}$ is strictly convex for $|t| < \epsilon$ and $\chi_{\epsilon} = |t|$ for $|t| \ge 0$. Then we define a smooth maximum by 
	\begin{equation*}
		\widetilde{\max}_{\epsilon}(x,y) := \frac{x + y + \chi_{\epsilon}(x-y)}{2}.
	\end{equation*}
	Note that $\widetilde{\max}_{\epsilon}(x,y) = \max(x,y)$ if $|x-y| \ge \epsilon$. Moreover, the smooth maximum of two $C^k$-smooth (strictly) convex functions is $C^k$-smooth (strictly) convex. Let $\dist(x,\partial \O) := \inf \{ \norm{x-y} : y \in \partial \O \}$. Let $\delta : \RR^n \rightarrow \RR$ be the distance function of $\O$ defined by
	\[
		\delta(x) := 
		\begin{cases}
			- \dist(x, \partial \O)  & \text{ if } x \in \O \\
			\dist(x, \partial \O) & \text{ if } x \in \O^{\complement}. 
		\end{cases}
	\]	
	Define $\O_{\epsilon} := \{ x \in \RR^n : \delta(x) < \epsilon \}$.

	\begin{thm} \label{strictly convex defining function}
		Let $\O \subset \subset \RR^n$ be a convex domain with $C^k(k \ge 1)$-smooth  boundary. Then there exists a $C^k$-smooth function $\rho : \RR^n \rightarrow \RR$ satisfying 
		\begin{itemize}
			\item[1.] $\rho$ is a defining function of $\O$.
			\item[2.] $\rho$ is strictly convex on $\RR^n \setminus \partial \O$.
		\end{itemize}
	\end{thm}

	The author would like to acknowledge that he learned the formulation of this theorem as well as the proof from Professor N. Shcherbina.
	
	\begin{proof}
		We may assume that $0 \in \O$. Consider the Minkowski function of $\O$. Define a function $f : \RR^n \rightarrow \RR$ defined by \\
		\[
		f(x) = 
			\begin{cases}
				\lambda, \text{ where } \frac{1}{\lambda}x \in \partial \O & \text{if } x \neq 0 \\
				0 & \text{if } x = 0.
			\end{cases}
		\]
		Then $f$ is well-defined,  $C^k$-smooth on $\RR^n \setminus \{0\}$, and $\nabla f \neq 0$ on $\partial \O$. Convexity of $\O$ implies that $f$ is a convex function. Consequently, $\sigma := f - 1$ is a convex defining function of $\O$. \\
		
		We first construct a defining function of $\O$ which is strictly convex inside of $\O$. For all $\epsilon > 0$ such that $0 \in \O_{-\epsilon}$, we claim that there exist $\delta_1, \delta_2 > 0$ such that $\delta_1 \norm{x}^2 - \delta_2 < 0$ on $\partial \O$, and $\delta_1 \norm{x}^2 - \delta_2 > \sigma$ on $\overline{\O}_{-\epsilon}$. 
		Let $S_{\epsilon} := \max \{ \sigma(x) : x \in \overline{\O}_{-\epsilon} \} < 0 $, 
		$m_{\epsilon} := \min \{ \norm{x}^2 : x \in \overline{\O}_{-\epsilon} \} > 0$,
		$M_{\epsilon} := \max \{ \norm{x}^2 : x \in \partial \O \} > 0$.
		Then one may choose sufficiently small $\delta_1, \delta_2 > 0$ such that 
		$$ \delta_1 M_{\epsilon} < \delta_2 < \delta_1 m_{\epsilon} + (-S_{\epsilon}). $$
		Therefore 
		\begin{equation*}
			\delta_1 \norm{x}^2 - \delta_2 \ge \delta_1 m_{\epsilon} - \delta_2 > S_{\epsilon} > \sigma(x)
		\end{equation*}
		on $\overline{\O}_{-\epsilon}$, and 
		\begin{equation*}
			\delta_1 \norm{x}^2 - \delta_2 \le \delta_1 M_{\epsilon} - \delta_2 < 0
		\end{equation*}		
		on $\partial \O$. The claim is proved.
	
		Let $V_{\epsilon} := \{ x \in \RR^n : \sigma - (\delta_1 \norm{x}^2 - \delta_2) < 0 \}$ and $m := \min\{ \dist(\partial V_{\epsilon}, \partial \O_{-\epsilon}),  \dist(\partial V_{\epsilon}, \partial \O) \}$. Define 
		$$\rho_{\epsilon} := \widetilde{\max}_{\frac{m}{2}}(\sigma, \delta_1 \norm{x}^2 - \delta_2).$$ 
		Then $\rho_{\epsilon}$ is a $C^k$-smooth function on $\RR^n$ and $\rho_{\epsilon} = \delta_1 \norm{x}^2 - \delta_2$ on $\O_{-\epsilon}$ which is strictly convex. Assume that $0 \in \O_{1}$ (the number 1 is not important here and one may choose a number smaller than 1). Let $U_{\epsilon} := \{ x \in \RR^n : |\delta(x)| < \epsilon \}$ and for $j \in \NN$
		$$ c_j := \sup_{\overline{U}_{j}} \sup_{0 \le |\alpha| \le k}  \left| \frac{\partial^{\alpha}}{\partial x^{\alpha}} \rho_{\frac{1}{j}}(x) \right|, \phantom{aa} \eta_j := \frac{1}{2^j c_j} $$
		Here, we used multi-index $\alpha = (\alpha_1, \cdots, \alpha_n) \in \NN^n_0$ and $|\alpha| = \alpha_1 + \cdots + \alpha_n$, $\frac{\partial^{\alpha}}{\partial x^{\alpha}} = \frac{\partial^{\alpha}}{\partial x_1^{\alpha_1} \cdots \partial x_n^{\alpha_n}}$.
		Define 
		$$ \sigma_1(x) := \sum_{j=1}^{\infty} \eta_j \rho_{\frac{1}{j}}(x). $$
		Then for all $j \in \NN$, $0 \le |\alpha| \le k$, since $\left| \frac{\partial^{\alpha}}{\partial x^{\alpha}} \rho_{\frac{1}{j}}(x) \right| \le \frac{1}{2^j}$ on $\overline{U}_j$ and $\sum_{j=1}^{\infty} \frac{1}{2^j}$ converges, 
		
		$$\sum_{j=1}^{\infty} \eta_j \frac{\partial^{\alpha}}{\partial x^{\alpha}} \rho_{\frac{1}{j}}(x)$$
		is uniformly convergent and hence continuous on $\RR^n$. Note that $\nabla \sigma_1 = \sum_{j=1}^{\infty} \eta_j \nabla \rho_{\frac{1}{j}} \neq 0$ on $\partial \O$. Therefore $\sigma_1$ is a well-defined $C^k$-smooth function on $\RR^n$ which is strictly convex on $\O$. \\

		Now, we consider the outside of $\O$. The argument is similar to that above. There exist $\widetilde{\delta}_1, \widetilde{\delta}_2 > 0$ such that $\sigma(x) + \widetilde{\delta}_1 \norm{x}^2 - \widetilde{\delta}_2 < 0$ on $\partial \O$, and $\sigma(x) + \widetilde{\delta}_1 \norm{x}^2 - \widetilde{\delta}_2 > 0$ on $\partial \O_{\epsilon}$. Let 
		$\widetilde{V}_{\epsilon} := \{ x \in \RR^n : \sigma + \widetilde{\delta}_1 \norm{x}^2 - \widetilde{\delta}_2 < 0 \}$ and 
		$\widetilde{m} := \min\{ \dist(\partial \widetilde{V}_{\epsilon}, \partial \O_{\epsilon}),  \dist(\partial \widetilde{V}_{\epsilon}, \partial \O) \}$. 
		Then 
		$$\widetilde{\rho}_{\epsilon} := \widetilde{\max}_{\frac{\widetilde{m}}{2}}(0, \sigma(x) + \widetilde{\delta}_1 \norm{x}^2 - \widetilde{\delta}_2)$$
		is $C^k$-smooth function on $\RR^n$ which is strictly convex on $\O_{\epsilon}^{\complement}$. Finally, for $\widetilde{\eta}_j > 0$ such that $\widetilde{\eta_j} \searrow 0$ sufficiently fast, 
		$$ \rho(x) := \sigma_1(x) + \sum_{j=1}^{\infty} \widetilde{\eta}_j \widetilde{\rho}_{\frac{1}{j}}(x) $$
		is the desired defining function of $\O$.
	\end{proof}

	\begin{cor} \label{convex domain index 1}
		For a convex domain $\O \subset \subset \CC^n$ with $C^k(k \ge 1)$-smooth boundary, $DF(\O)=1$ and $S(\O)=1$.
	\end{cor}
	\begin{proof}
		By Theorem \ref{strictly convex defining function}, there exists a defining function $\rho$ of $\O$ which is strictly convex on $\CC^n \setminus \partial \O$. Hence $\rho$ is strictly plurisubharmonic on $\CC^n \setminus \partial \O$. This implies that $DF(\O)=1$ and $S(\O)=1$.
	\end{proof}
	
	\begin{rmk}
		If the boundary regularity of a convex domain is $C^{\infty}$-smooth, then there is another way to prove Corollary \ref{convex domain index 1}. $\sigma$ in the proof of Theorem \ref{strictly convex defining function} is $C^{\infty}$-smooth plurisubharmonic defining function of $\O$. Hence the theorem by Forn{\ae}ss and Herbig (\cite{ForHer2008}) implies $DF(\O)=1$ and $S(\O)=1$.
	\end{rmk}

	\vspace{5mm}

	
	\textit{\bf Acknowledgements}: The author would like to express his deep gratitude to Professor Kang-Tae Kim for valuable guidance and encouragements, and to Professor N. Shcherbina for fruitful conversations.
	

	\vspace{5mm}

	\bibliographystyle{amsplain}

\end{document}